\documentclass{arxiv}
\usepackage{amsmath, amssymb, amsthm, graphicx}
\usepackage{comment}
\newtheorem{theorem}{Theorem}[section]

\theoremstyle{definition}

\newtheorem{criterion}[]{Criterion}
\newtheoremstyle{named}{}{}{\itshape}{}{\bfseries}{.}{.5em}{\thmnote{#3's }#1} \theoremstyle{named} 

\theoremstyle{remark}
\newtheorem{remark}[theorem]{Remark}

\numberwithin{equation}{section}

\numberwithin{equation}{section}

\title{Polygonal homographic orbits in spaces of constant curvature}

\author{Pieter Tibboel}
\address{\em Department of Mathematics and Statistics \\Chongqing University, China}
\email{\em Pieter.Tibboel@gmail.com}
\begin{document}
\maketitle
\begin{abstract}
  We prove that the geometry of the 2-dimensional $n$-body problem for spaces of constant curvature $\kappa\neq 0$, $n\geq 3$, does not allow for polygonal homographic solutions, provided that the corresponding orbits are irregular polygons of non-constant size.
\end{abstract}

\section{Introduction}
By the $n$-body problem for spaces of constant curvature, we mean the problem of describing the dynamics of $n$ point particles in a space of constant curvature $\kappa\neq 0$. Polygonal homographic solutions are solutions to such a problem for which the point particles describe the vertices of a polygon that retains its shape over time. \\
This paper is inspired by work done by Diacu \cite{D2}, but research on this type of problem goes back as far as the 1830s when Bolyai \cite{BB} and Lobachevsky \cite{Lo} independently proposed a curved 2-body problem in hyperbolic space $\mathbb{H}^{3}$. Since then, the problem has been studied by outstanding mathematicians such as Dirichlet, Schering \cite{S1}, \cite{S2}, Killing \cite{K1}, \cite{K2}, \cite{K3} and Liebmann \cite{L1}, \cite{L2}, \cite{L3}. More recent results were obtained by Cari\~nena, Ra\~nada, Santander \cite{CRS}, Diacu \cite{D1}, \cite{D2}, Diacu, P\'erez-Chavela \cite{DP} and Diacu, P\'erez-Chavela, Santoprete \cite{DPS1}, \cite{DPS2}, \cite{DPS3}. For a more detailed historical overview, please see \cite{D2} or \cite{DPS1}.\\
This paper is about the existence of polygonal homographic solutions. An important reason for studying these objects is that they may give us information about the geometry of our universe. For example, Diacu, P\'erez-Chavela, and M. Santoprete showed in \cite{DPS1}, \cite{DPS2} that homographic equilateral triangles of non-equal masses (the so-called Lagrangian solutions) exist only in Euclidean space and that in the curved 3-body problem the masses must be equal.
Because our universe contains bodies of non-equal masses that traverse homographic equilateral triangles, we may assume that, at least locally, our universe is Euclidean (see \cite{D2}, \cite{DPS1}, \cite{DPS2}). \\
As regular homographic orbits for spaces of constant curvature have been well investigated (see \cite{D2}, \cite{DPS1}, \cite{DPS2}, \cite{DPS3}), we will focus in this paper on irregular homographic orbits. The most general theorem sofar regarding the existence of irregular polygonal homographic orbits is by Diacu, who proved in \cite{D2} for $n=3$:
\begin{theorem}\label{NonEquilateralTriangles}
  Consider the curved $3$-body problem, given by equations~(\ref{EquationsOfMotion})
  with $n=3$ and masses $m_{1}$, $m_{2}$, $m_{3}>0$.
  These equations admit no homographic orbits given by scalene non-equilateral triangles for $\kappa<0$. For $\kappa>0$, they do not admit such solutions either if the bodies stay away from the equator.
\end{theorem}
What we will prove in this paper is the new result:
\begin{theorem}\label{NonExistenceOfIrregularPolygons}
  Consider the curved $n$-body problem, given by equations (\ref{EquationsOfMotion}),
  with $n\geq 3$ and masses $m_{1}$,...,$m_{n}>0$.
  These equations do not allow for homographic orbits given by irregular polygons if $z$ is not constant,
  where $z$ is the same as in (\ref{Representation}).
\end{theorem}
\begin{remark}
  Of particular interest is that our proof relies on the irregularity of the polygons alone and does not make use of the values of the masses. This means that the non-existence of the orbits is determined solely by the geometry of the space.
\end{remark}
Before we can prove Theorem~\ref{NonExistenceOfIrregularPolygons}, we need to formulate a criterion due to Diacu \cite{D2}, which will play a key role in our proof. This will be done in the next section.
\section{Diacu's Criterion}
In this section, we will formulate a criterion that gives necessary and sufficient conditions for the existence of polygonal homographic orbits. The notation used in this paper was introduced by Diacu in \cite{D1} and makes it possible to use a unified formulation for both the positive and negative constant curvature case, thus greatly simplifying calculations.\\
Consider the curved $n$-body problem of $n$ point particles. \\We
will denote their masses to be $m_{1}$, $m_{2}$,..., $m_{n}>0$ and
their positions by the vectors
$\textbf{q}_{i}=(x_{i},y_{i},z_{i})\in\textbf{M}_{\kappa}^{2}$,
$i=\overline{1,n}$, \\where
$\textbf{M}_{\kappa}^{2}=\{(x,y,z)\in\mathbb{R}^{3}|\textrm{
}\kappa(x^{2}+y^{2}+\sigma z^{2})=1\}$ and
\begin{align*}
    \sigma=\begin{cases}
      \hspace{0.25cm}1 &\textrm{ for }\kappa\geq 0\\
      -1 &\textrm{ for }\kappa<0
    \end{cases}
  \end{align*}
Furthermore, consider for 3-dimensional vectors $\textbf{a}=(a_{x},a_{y},a_{z})$, $\textbf{b}=(b_{x},b_{y},b_{z})$ the inner product
  \begin{align*}
    \textbf{a}\odot\textbf{b}=a_{x}b_{x}+a_{y}b_{y}+\sigma
    a_{z}b_{z}.
  \end{align*}
Then, following \cite{D1}, \cite{D2}, \cite{DPS1}, \cite{DPS2} and \cite{DPS3}, we define the equations of motion for the $n$-body problem as the dynamical system described by
\begin{align}\label{EquationsOfMotion}
  \ddot{\textbf{q}}_{i}=\sum\limits_{j=1, j\neq i}^{n}\frac{m_{j}|\kappa|^{\frac{3}{2}}[\textbf{q}_{j}-(\kappa\textbf{q}_{i}\odot\textbf{q}_{j})\textbf{q}_{i}]}{[\sigma-\sigma(\kappa\textbf{q}_{i}\odot\textbf{q}_{j})^{2}]^{\frac{3}{2}}}-(\kappa\dot{\textbf{q}}_{i}\odot\dot{\textbf{q}}_{i})\textbf{q}_{i},\textrm{ }i=\overline{1,n}
\end{align}
We call the solution of (\ref{EquationsOfMotion}) a \textit{polygonal homographic solution}, or a \textit{polygonal homographic orbit}, in accordance with  \cite{D2}, if it can be represented as
\begin{align}\label{Representation}
  \textbf{q}=(\textbf{q}_{1},...,\textbf{q}_{n}),\textrm{ }\textbf{q}_{i}=(x_{i},y_{i},z_{i})\hspace{2cm}\\
  x_{i}=r\cos{(\omega+\alpha_{i})},\textrm{ }y_{i}=r\sin{(\omega+\alpha_{i})},\textrm{ }z_{i}=z,\textrm{ }\overline{1,n},\nonumber
\end{align}
where $0\leq\alpha_{1}<\alpha_{2}<...<\alpha_{n}<2\pi$ are constants, the function $z=z(t)$ satisfies $z^{2}=\sigma\kappa^{-1}-\sigma r^{2}$, $r:=r(t)$ is the \textit{size function} and $\omega:=\omega(t)$ is the \textit{angular function}. \\
The main tool of our proof of Theorem~\ref{NonExistenceOfIrregularPolygons} is Criterion 1 as formulated in \cite{D2}, which states:
\begin{criterion}[Diacu's Criterion] Consider $n\geq 3$ bodies of masses $m_{1}$, $m_{2}$,...,$m_{n}>0$ moving on the surface $\textbf{M}_{\kappa}^{2}$. The necessary and sufficient conditions for a polygonal homographic orbit as described in (\ref{Representation}) to be a solution of equations (\ref{EquationsOfMotion}) are given by the equations
\begin{align}
  \delta_{1}=\delta_{2}=...=\delta_{n}\textrm{ and }\gamma_{1}=\gamma_{2}=...=\gamma_{n},
\end{align}
where
\begin{align}
  \delta_{i}=\sum\limits_{j=1,j\neq i}^{n}m_{j}\mu_{ji},\textrm{ }\gamma_{i}=\sum\limits_{j=1,j\neq i}^{n}m_{j}\nu_{ji},\textrm{ }, i=\overline{1,n}, \\
  \mu_{ji}=\frac{1}{c_{ji}^{\frac{1}{2}}(2-c_{ji}\kappa r^{2})^{\frac{3}{2}}},\textrm{ }\nu_{ji}=\frac{s_{ji}}{c_{ji}^{\frac{3}{2}}(2-c_{ji}\kappa r^{2})^{\frac{3}{2}}}\label{DefinitionmuANDnu}\\
  s_{ji}=\sin{(\alpha_{j}-\alpha_{i})},\textrm{ }c_{ji}=1-\cos{(\alpha_{j}-\alpha_{i})},\textrm{ }i,j=\overline{1,n},\textrm{ }, i\neq j.\label{DefinitionsANDc}
\end{align}
\end{criterion}
\section{Proof of Theorem~\ref{NonExistenceOfIrregularPolygons}}
With Diacu's Criterion in place, we can now prove Theorem~\ref{NonExistenceOfIrregularPolygons}:
\begin{proof}
  For the purpose of convenience later on in the proof, we define \\$\rho:=\kappa r^{2}$ and $\alpha_{n+1}:=\alpha_{1}+2\pi$. Furthermore, we choose $\alpha_{2}-\alpha_{1}\leq\alpha_{i+1}-\alpha_{i}$, $\alpha_{i}$ as in (\ref{Representation}), $i=\overline{1,n}$. This can be done by choosing a suitable $x$-axis and $y$-axis. \\
  Assume that (\ref{EquationsOfMotion}) does allow for homographic orbits given by an irregular polygon if $z$ is not constant. Then according to Criterion 1, $\delta_{1}-\delta_{2}=0$ and $\gamma_{1}-\gamma_{2}=0$, so
  \begin{align}\label{FirstEqualities1}
    0=\delta_{1}-\delta_{2}=(m_{2}-m_{1})\mu_{21}+\sum\limits_{j=3}^{n}m_{j}(\mu_{j1}-\mu_{j2}),
  \end{align}
  as $\mu_{21}=\mu_{12}$ and
  \begin{align}\label{FirstEqualities2}
    0=\gamma_{1}-\gamma_{2}=(m_{2}+m_{1})\nu_{21}+\sum\limits_{j=3}^{n}m_{j}(\nu_{j1}-\nu_{j2}),
  \end{align}
  as $\nu_{21}=-\nu_{12}$. \\
  Note that, using (\ref{DefinitionmuANDnu}) and (\ref{DefinitionsANDc}), $\nu_{ji}=\frac{s_{ji}}{c_{ji}}\mu_{ji}$, so (\ref{FirstEqualities2}) can be rewritten as
    \begin{align}\label{FirstEqualities2a}
        0=(m_{2}+m_{1})\frac{s_{21}}{c_{21}}\mu_{21}+\sum\limits_{j=3}^{n}m_{j}\left(\frac{s_{j1}}{c_{j1}}\mu_{j1}-\frac{s_{j2}}{c_{j2}}\mu_{j2}\right).
  \end{align}
  Using the definition of $\rho$ and (\ref{DefinitionmuANDnu}), we write
  \begin{align}\label{Justmu}
    \mu_{ji}(\rho)=\frac{1}{c_{ji}^{\frac{1}{2}}(2-c_{ji}\rho)^{\frac{3}{2}}}
  \end{align}
  As $z$ is not constant, $\rho$ is not constant either and we may
  take the $k$th derivative of (\ref{Justmu}), $k\in\mathbb{N}\cup\{0\}$, with respect to
  $\rho$ to obtain
  \begin{align}\label{Derivativemu}
    \mu_{ji}^{(k)}(\rho)=\left(\prod_{l=0}^{k-1}\left(\frac{3}{2}+l\right)\right)\frac{c_{ji}^{\frac{1}{2}+k}}{(2-c_{ji}\rho)^{\frac{3}{2}+k}}
  \end{align}
  Taking the $k$th derivatives with respect to $\rho$ of (\ref{FirstEqualities1}) and (\ref{FirstEqualities2a}), inserting (\ref{Derivativemu}) and dividing everything by $\left(\prod_{l=0}^{k-1}\left(\frac{3}{2}+l\right)\right)$ then gives
  \begin{align}\label{SecondEqualities1}
    0=(m_{1}-m_{2})\frac{c_{21}^{\frac{1}{2}+k}}{(2-c_{21}\rho)^{\frac{3}{2}+k}}+\sum\limits_{j=3}^{n}m_{j}\left(\frac{c_{j1}^{\frac{1}{2}+k}}{(2-c_{j1}\rho)^{\frac{3}{2}+k}}-\frac{c_{j2}^{\frac{1}{2}+k}}{(2-c_{j2}\rho)^{\frac{3}{2}+k}}\right),
  \end{align}
  and
  \begin{align}\label{SecondEqualities2}
    \hspace{-0.5cm}0=(m_{1}+m_{2})\frac{s_{21}}{c_{21}}\frac{c_{21}^{\frac{1}{2}+k}}{(2-c_{21}\rho)^{\frac{3}{2}+k}}+\sum\limits_{j=1}^{n}m_{j}\left(\frac{s_{j1}}{c_{j1}}\frac{c_{j1}^{\frac{1}{2}+k}}{(2-c_{j1}\rho)^{\frac{3}{2}+k}}-\frac{s_{j2}}{c_{j2}}\frac{c_{j2}^{\frac{1}{2}+k}}{(2-c_{j2}\rho)^{\frac{3}{2}+k}}\right).
  \end{align}
  We now fix $\rho$ and rewrite (\ref{SecondEqualities1}) and (\ref{SecondEqualities2}) defining
  \begin{align*}
    \frac{c_{ji}^{\frac{1}{2}+k}}{(2-c_{ji}\rho)^{\frac{3}{2}+k}}=a_{ji}(\rho)g_{ji}(\rho)^{k},
  \end{align*}
  where
  \begin{align}\label{DefinitionaANDg}
    a_{ji}(\rho)=\frac{c_{ji}^{\frac{1}{2}}}{(2-c_{ji}\rho)^{\frac{3}{2}}}\textrm{ and }g_{ji}(\rho)=\frac{c_{ji}}{(2-c_{ji}\rho)}
  \end{align}
  Then (\ref{SecondEqualities1}) and (\ref{SecondEqualities2}) become
  \begin{align}\label{SecondEqualities1I}
    0=(m_{1}-m_{2})a_{21}g_{21}^{k}+\sum\limits_{j=3}^{n}m_{j}\left(a_{j1}g_{j1}^{k}-a_{j2}g_{j2}^{k}\right),
  \end{align}
  and
  \begin{align}\label{SecondEqualities2I}
    0=(m_{1}+m_{2})\frac{s_{21}}{c_{21}}a_{21}g_{21}^{k}+\sum\limits_{j=3}^{n}m_{j}\left(\frac{s_{j1}}{c_{j1}}a_{j1}g_{j1}^{k}-\frac{s_{j2}}{c_{j2}}a_{j2}g_{j2}^{k}\right),
  \end{align}
  Note that the right hand sides of (\ref{SecondEqualities1I}) and (\ref{SecondEqualities2I}) are linear combinations of exponential functions in $k$, where $k$ can vary in $\mathbb{N}\cup\{0\}$. These exponential functions are linearly independent, provided their bases are distinct.
  The main idea of the proof will be to show that there is a function $g_{j1}^{k}$, $j\in\{1,...n\}$, in the linear combinations of (\ref{SecondEqualities1I}) and (\ref{SecondEqualities2I}) that does not cancel out in at least one of the equations (\ref{SecondEqualities1I}) and (\ref{SecondEqualities2I}). \\
  Let us assume the contrary.
  For a function $g_{j1}^{k}$ to be canceled out, we need to be able to construct a linear combination of $g_{j1}^{k}$ with other exponential functions in the linear combinations of (\ref{SecondEqualities1I}) and (\ref{SecondEqualities2I}), with identical bases. Such functions can be represented as $g_{u2}^{k}$, $g_{v1}^{k}$, $u$, $v\in\{1,...,n\}$ for which $g_{j1}=g_{u2}$, $g_{j1}=g_{v1}$.
   However, due to (\ref{DefinitionaANDg}), $g_{j1}=g_{u2}$ can only be the case if $c_{j1}=c_{u2}$ and $g_{j1}=g_{v1}$ can only be the case if $c_{j1}=c_{v1}$. \\
  We need to take a closer look at the case that $c_{j1}=c_{u2}$. After that, we will make some further remarks about $u$, $v\in\{1,...,n\}$ for which $c_{j1}=c_{u2}$ and $c_{j1}=c_{v1}$ hold
  respectively, which will then allow us to move on to the final part of our proof.\\
  For the case that $c_{j1}=c_{u2}$, we have that by (\ref{DefinitionsANDc}),
  \begin{align}\label{Possibility1}
    \alpha_{j}-\alpha_{1}=\alpha_{u}-\alpha_{2} (\textrm{ mod }2\pi),
  \end{align}
  or
  \begin{align}\label{Possibility2}
    \alpha_{j}-\alpha_{1}=\alpha_{2}-\alpha_{u} (\textrm{ mod }2\pi)
  \end{align}
  Consequently, this implies for (\ref{Possibility1}) that
  \begin{align}\label{Possibility1b}
    \alpha_{u}=\alpha_{j}+\alpha_{2}-\alpha_{1}>\alpha_{j}
  \end{align}
  Note that by construction
  \begin{align}\label{constructioninequality}
    \alpha_{2}-\alpha_{1}\leq\alpha_{j+1}-\alpha_{j}.
  \end{align}
  Combining (\ref{constructioninequality}) with (\ref{Possibility1b}), we get
  \begin{align}\label{Possibility1c}
    \alpha_{j}<\alpha_{u}=\alpha_{j}+\alpha_{2}-\alpha_{1}\leq\alpha_{j}+\alpha_{j+1}-\alpha_{j}=\alpha_{j+1}
  \end{align}
  So if (\ref{Possibility1}) holds, then $\alpha_{j}<\alpha_{u}\leq\alpha_{j+1}$, so $\alpha_{u}=\alpha_{j+1}$ and thus
  \\$\alpha_{j+1}-\alpha_{j}=\alpha_{2}-\alpha_{1}$.
  If we assume that for every $j$ a $u$ exists such that (\ref{Possibility1}) is true, then
  \begin{align}
    2\pi=\sum\limits_{j=1}^{n}(\alpha_{j+1}-\alpha_{j})=n(\alpha_{2}-\alpha_{1})=n(\alpha_{j+1}-\alpha_{j})
  \end{align}
  So $\alpha_{j+1}-\alpha_{j}=\frac{2\pi}{n}$ for all $j\in\{1,...,n\}$. However, we assumed that our polygon was not regular, so there must be at least one $j\in\{1,...,n\}$ for which (\ref{Possibility2}) holds and (\ref{Possibility1}) does not.
  This $j$ will give us the desired contradiction. \\
  Before moving on to the last part of the proof, we should make the following observations:
  If there is a $v$ such that $c_{j1}=c_{v1}$,
   then either $\alpha_{j}+\alpha_{v}=2\alpha_{1}$, or $\alpha_{j}+\alpha_{v}=2\alpha_{1}+2\pi$.
   Since $\alpha_{1}<\alpha_{j}$ and $\alpha_{1}<\alpha_{v}$, $\alpha_{j}+\alpha_{v}=2\alpha_{1}$
   would mean that $\alpha_{j}+\alpha_{v}>\alpha_{j}+\alpha_{v}$. \\
   If $\alpha_{j}-\alpha_{1}=\alpha_{2}-\alpha_{u}(\textrm{ mod }2\pi)$, then $\alpha_{j}+\alpha_{u}=\alpha_{1}+\alpha_{2}(\textrm{ mod }2\pi)$. \\If $\alpha_{j}+\alpha_{u}=\alpha_{1}+\alpha_{2}$, then, as $\alpha_{u}>\alpha_{j}$, $\alpha_{u}=\alpha_{2}$ and $\alpha_{j}=\alpha_{1}$, as $\alpha_{1}$ and $\alpha_{2}$ are the smallest angles available. However, $\alpha_{j}\neq\alpha_{1}$, so $\alpha_{u}=-\alpha_{j}+\alpha_{1}+\alpha_{2}+2\pi$. \\
   From this we derive:
  \begin{itemize}
    \item[\textit{I}.] If there is a $v$ such that $c_{j1}=c_{v1}$, then $v$ is
    unique, $\alpha_{j}-\alpha_{1}=\alpha_{1}-\alpha_{v}(\textrm{ mod }2\pi)$ and $s_{j1}=-s_{v1}$.
    \item[\textit{II}.] If there is a $u$ such that
    (\ref{Possibility2}) holds, then $u$ is
    unique, $\alpha_{j}-\alpha_{1}=\alpha_{2}-\alpha_{u}(\textrm{ mod }2\pi)$ and $s_{j1}=-s_{u2}$.
    \item[\textit{III}.] $c_{st}=c_{\tilde{s}\tilde{t}}$ implies that $a_{st}=a_{\tilde{s}\tilde{t}}$, $s$, $t$, $\tilde{s}$, $\tilde{t}\in\{1,...,n\}$, $s\neq t$, $\tilde{s}\neq \tilde{t}$.
  \end{itemize}
  For the aforementioned value $j$ there are now three possibilities:
  \begin{itemize}
    \item[1.]  For $j$ there exists neither a value $u$ such that (\ref{Possibility2}) holds nor is there a value $v$ such that $c_{j1}=c_{v1}$.
    \item[2.] For $j$ there exists either a unique value $u$ such that (\ref{Possibility2}) holds or there is a unique value $v$ such that $c_{j1}=c_{v1}$.
    \item[3.] For $j$ there exists both a unique value $u$ such that (\ref{Possibility2}) holds and a unique value $v$ such that $c_{j1}=c_{v1}$.
  \end{itemize}
  For the first possibility, $g_{j1}^{k}$ does not cancel out, as there are no other terms in either linear combination (\ref{SecondEqualities1I}) or (\ref{SecondEqualities2I}) that depend on it, which means we have a contradiction. \\
  Let us consider the second possibility. If there exists a unique value $u$ such that (\ref{Possibility1}) holds, using items \textit{II} and \textit{III} in the list above, in (\ref{SecondEqualities2I}) the coefficients of $g_{j1}^{k}$ and $g_{u2}^{k}$
  add up to,
  \begin{align}\label{aformula}
    a_{j1}\left(\frac{s_{j1}}{c_{j1}}m_{j}-\frac{s_{u2}}{c_{u2}}m_{u}\right)=a_{j1}\left(\frac{s_{j1}}{c_{j1}}m_{j}--\frac{s_{j1}}{c_{j1}}m_{u}\right)=a_{j1}\frac{s_{j1}}{c_{j1}}(m_{j}+m_{u}).
  \end{align}
  We should check that $s_{j1}$ is not equal to zero in (\ref{aformula}): \\
  $s_{j1}=0$ if and only if $\sin{(\alpha_{j}-\alpha_{1})}=0$, if and only if $\alpha_{j}-\alpha_{1}=0 (\textrm{ mod }\pi)$. \\
  Because $\alpha_{j}\neq\alpha _{1}$ by the definition of $\mu_{j1}$ and $\nu_{j1}$, the only other possibility is that  $\alpha_{j}-\alpha_{1}=\pi$, which would mean that $\pi=\alpha_{2}-\alpha_{u} (\textrm{ mod }2\pi)$. However, then $\alpha_{u}-\alpha_{2}=\pi (\textrm{ mod }2\pi)$, which implies that (\ref{Possibility1}) holds for $j$ and $u$, which we chose not to be the case. \\
  Thus $s_{j1}\neq 0$, which means that we can deduce the desired
  contradiction from equation (\ref{SecondEqualities2I}). \\
  Let us look at the second case when there exists a unique value $v$ such that $c_{j1}=c_{v1}$. Using item~\textit{III} of the list above, in (\ref{SecondEqualities1I}), the coefficients of $g_{j1}^{k}$ and $g_{v1}^{k}$ add up to $a_{j1}(m_{j}+m_{v})$, which gives the desired contradiction. \\
  For possibility 3, using items \textit{I}, \textit{II} and \textit{III} of the list above, the coefficients of $g_{j1}^{k}$, $g_{u2}^{k}$ and $g_{v1}^{k}$ add up to $a_{j1}(m_{j}+m_{v}-m_{u})$
  in (\ref{SecondEqualities1I}) and they add up to \\$a_{j1}\frac{s_{j1}}{c_{j1}}(m_{j}-m_{v}+m_{u})$
  in (\ref{SecondEqualities2I}). As $m_{j}$, $m_{v}$ and $m_{u}$ are all positive, the coefficients of $g_{j1}^{k}$, $g_{v1}^{k}$ and $g_{u2}^{k}$ do not cancel out in at least one of the equations (\ref{SecondEqualities1I}) and (\ref{SecondEqualities2I}),
  which gives us the desired contradiction.
\end{proof}

\end{document}